\documentclass[reqno]{alea2}
\usepackage{natbib}
\usepackage{fancyhdr}
\usepackage{graphicx}

\usepackage{amssymb}
\usepackage{amsmath}

\renewcommand{\cite}{\citet}

\makeatletter \@addtoreset{equation}{section} \makeatother

\renewcommand\theequation{\thesection.\arabic{equation}}
\renewcommand\thefigure{\thesection.\@arabic\c@figure}
\renewcommand\thetable{\thesection.\@arabic\c@table}

\theoremstyle{plain}
\newtheorem{theorem}{Theorem}[section]                                          \newtheorem{proposition}[theorem]{Proposition}
\newtheorem{lemma}[theorem]{Lemma}

\newtheorem{Cl}{Claim}[section]
\theoremstyle{definition}

\theoremstyle{remark}

\newcommand{\N}{\mathbb{N}}
\newcommand{\R}{\mathbb{R}}
\newcommand{\T}{\mathbb{P}}
\newcommand{\Z}{\mathbb{Z}}
\newcommand{\C}{\mathbb{C}}

\newcommand{\bal}{\begin{align}}
\newcommand{\eal}{\end{align}}
\newcommand{\beq}{\begin{equation}}
\newcommand{\eeq}{\end{equation}}
\newcommand{\bpf}{\begin{proof}}
\newcommand{\epf}{\end{proof}}
\newcommand{\bsp}{\begin{split}}
\newcommand{\esp}{\end{split}}
\newcommand{\bg}{\begin{gathered}}
\newcommand{\eg}{\end{gathered}}

\begin{document}

\date{, accepted }

\author{Cristina Tone}
\address{Department of Mathematics, University of Louisville, 328 NS, Louisville, Kentucky 40292}
\email{cristina.tone@louisville.edu}

\title[CLTs for Hilbert-space valued random fields under $\rho'$-mixing]{Central limit theorems for Hilbert-space valued random fields satisfying a strong mixing condition}

\begin{abstract}
In this paper we study the asymptotic normality of the normalized
partial sum of a Hilbert-space valued strictly stationary random field satisfying the interlaced $\rho'$-mixing condition.
\end{abstract}

\maketitle

\section{Introduction}\label{s:intro}

In the literature about Hilbert-valued random sequences under mixing conditions, progress has been made by \citet{Mal},  \citet{FM}, and \citet{MPU}.
\citet{DM} established a central limit theorem 
and its weak invariance principle for Hilbert-valued strictly stationary sequences
under a ''projective criterion.'' In this way, they recovered the special case
of Hilbert-valued martingale difference sequences, and under a strong mixing
condition involving the whole past of the process and just one ''future''
observation at a time, they gave the nonergodic version of the result of
\citet{MPU}. Later on,  \citet{FM} proved a 
central limit theorem for a Hilbert-space valued strictly stationary,
strongly mixing sequence, where the mixing coefficients involve the whole past of the
process and just two ''future'' observations at a time, by using the Bernstein blocking technique and
approximations by martingale differences.

This paper will present a central limit theorem for strictly
stationary Hilbert-space valued random fields satisfying the $\rho'$-mixing condition.
We proceed by proving in Theorem \ref{r} a central limit theorem for a $\rho'$-mixing
strictly stationary random field of real-valued random variables, by the use of the Bernstein blocking technique. Next, in Theorem \ref{rm} we extend the
real-valued case to a random field of $m$-dimensional random vectors,
$m\geq1$, satisfying the same mixing condition. Finally, being able to prove the tightness condition in Theorem \ref{hilb}, we extend the
finite-dimensional case even further to a (infinite-dimensional) Hilbert
space-valued strictly stationary random field in the presence of the $\rho'$-mixing condition.

\section{Preliminary Material}\label{s:pre}
For the clarity of the proofs of the three theorems mentioned above, relevant definitions,
notations and basic background information will be given first.

Let $(\Omega, \mathcal{F}, \T)$ be a probability space.
Suppose $H$ is a separable real Hilbert space with inner product
$\langle \cdot    ,\cdot \rangle$ and norm $\|\cdot \|_H$.
Let $\mathcal{H}$ be the $\sigma$-field generated by the class of
all open subsets of $H$.
Let $\{ e_k \}_{k\geq 1}$ be an orthonormal basis for the Hilbert
space $H$.
Then for every $x \in H $, we denote by $x_k$ the kth
coordinate of $x$, defined by $ x_k= \langle x ,e_k
\rangle$, $k\geq 1.$
Also, for every $x \in H $ and every
$N\geq 1$ we set \[r_N^2(x)=\sum_{k=N}^{\infty} x_k^2
=\sum_{k=N}^{\infty} \langle x,e_k \rangle^2.\]
For any given $H$-valued random variable X with $EX=0_H$
and $E\|X \|_H^2 <\infty$, represent X by
\[X=\sum_{k=1}^{\infty}X_k e_k,\]
where $X_1, X_2, X_3, \ldots$
are real-valued random variables having $EX_k=0$ and $EX_k^2< \infty,
\ \forall \, k\geq 1$ (in fact, $\sum_{k=1}^{\infty}EX_k^2= E\|X
\|_H^2 <\infty$).
Then the ``covariance operator" (defined relative to the given orthonormal basis) for the (centered) $H$-valued random variable X
can be thought of as represented by the $\N\times \N$ ``covariance
matrix" $\Sigma :=(\sigma_{ij}, i\geq1, j\geq 1)$,
where $ \sigma_{ij}:=EX_iX_j. $
\begin{lemma}\label{hl1}
Let $\mathcal{P}_0$ be a class of probability measures on $(H,
\mathcal{H})$ satisfying the following conditions: \[\sup_{P \in
\mathcal{P}_0} \int_{H}  r_1^2(x)dP(x)<\infty, \text{ and}\]
\begin{equation*} \lim_{N\rightarrow \infty} \sup_{P \in
\mathcal{P}_0} \int_{H}  r_N^2(x)dP(x)=0. \end{equation*}
Then $\mathcal{P}_0$ is tight.
\end{lemma}
\noindent For the proof of the lemma, see \citet{Laha}, Theorem 7.5.1.\\

\noindent For any two $\sigma$-fields $\mathcal{A}, \ \mathcal{B} \subseteq
\mathcal{F}$, define now the strong mixing coefficient
\begin{equation*}
\alpha(\mathcal{A}, \mathcal{B}):=\sup_{A \in
\mathcal{A}, B \in \mathcal{B}}|P(A\cap B)-P(A)P(B)|, \end{equation*} and the
maximal coefficient of correlation
\begin{equation*}
\rho(\mathcal{A}, \mathcal{B}):=\sup|Corr(f, g)|, \text{ } f \in
L^2_{\text{real}}(\mathcal{A}), \text{ } g \in
L^2_{\text{real}}(\mathcal{B}).\end{equation*}
Suppose $d$ is a positive integer and $X:=(X_k, k\in
\Z^d)$ is a strictly stationary random field.
In this context, for each positive integer $n$, define the following
quantity:
\begin{equation*}\label{s1.8}
\alpha(n):=\alpha(X, n):=\sup \alpha(\sigma(X_k, k\in Q), \sigma(X_k, k \in S)),\end{equation*}
where the supremum is taken over all pairs of
nonempty,
disjoint sets $ Q$, $S \subset \Z^d$ with the
following property: There exist $u \in \{1, 2, \ldots, d\}$ and $j \in \Z $ such that
$Q \subset \{k:=(k_1, k_2, \ldots, k_d) \in \Z^d: k_u \leq j \}$ and
$S \subset \{k:=(k_1, k_2, \ldots, k_d) \in \Z^d: k_u \geq j+n \}$.

The random field $X:=(X_k, k\in \Z^d)$ is said to be ``strongly
mixing" (or ``$\alpha$-mixing") if $\alpha(n)\rightarrow 0$ as
$n\rightarrow\infty$.

Also, for each positive integer $n$, define the following quantity:
\begin{equation*}\label{s1.9} \rho'(n):=\rho'(X, n):=\sup \rho(\sigma(X_k, k\in Q), \sigma(X_k, k \in S)),\end{equation*} where the supremum is taken over all pairs of nonempty,
finite disjoint sets $Q$, $S \subset \Z^d$ with the
following property: There exist $u \in \{1, 2, \ldots, d\}$ and nonempty disjoint sets $A$, $B \subset \Z$, with $dist(A, B):=\min_{a\in A, b\in B}|a-b|\geq n$ such that $Q \subset \{k:=(k_1, k_2, \ldots, k_d) \in \Z^d: k_u \in A \}$ and $S \subset \{k:=(k_1, k_2, \ldots, k_d) \in \Z^d: k_u \in B \}$.

The random field $X:=(X_k, k\in \Z^d)$ is said to be
``$\rho'$-mixing" if $\rho'(n)\rightarrow 0$ as
$n\rightarrow\infty$.

Again, suppose $d$ is a positive integer, and
suppose $X:=(X_k, k\in \Z^d)$ is a strictly stationary Hilbert-space random field. Elements of $\N^d$ will be denoted by $L:=(L_1, L_2, \ldots, L_d)$. For any $L\in \N^d$, define the ``rectangular sum":
\begin{equation*}\label{s1.11}
S_{L}=S(X, L):=\sum_{k}{X_k},\end{equation*}where the sum is taken
over all $d$-tuples $k:=(k_1, k_2, \ldots, k_d) \in \N^d$ such that $1 \leq k_u \leq L_u \text{ for all } u\in \{1, 2, \ldots, d \}$. Thus $S(X, L)$ is the sum of $L_1\cdot L_2\cdot \ldots \cdot L_d$ of the $X_k's$.

\begin{proposition}\label{app1} Suppose $d$ is a positive integer.

 (I)  Suppose $(a(k), \ k \in \N^d)$ is an array of real (or complex)
numbers and $b$ is a real (or complex) number. Suppose that for
every $u \in \{1, 2, \ldots, d\}$ and every sequence $\left(
L^{(n)}, \ n \in \N \right)$ of elements of $\N^d$ such that
$L^{(n)}_u=n $ for all $n\geq 1$, and $L^{(n)}_v \rightarrow \infty$
as $n\rightarrow \infty$, $\forall \,$ $v \in \{1, 2, \ldots,
d\}\setminus\{u\}$, one has that $\lim_{n\rightarrow \infty}a \left( L^{(n)}
\right)=b.$ Then $a(L)\rightarrow b$ as $\min\{L_1, L_2, \ldots,
L_d\}\rightarrow\infty$.

(II)  Suppose $(\mu(k), \ k \in \N^d)$ is an array of
probability measures on $(S, \mathcal{S})$, where $(S, d)$ is a complete separable metric space and $\mathcal{S}$ is the $\sigma$-field on $S$ generated by the open balls in $S$ in the given metric $d$. Suppose $\nu$ is a
probability measure on $(S, \mathcal{S})$ and that for every $u
\in \{1, 2, \ldots, d\}$ and every sequence $(L^{(n)}, \text{ } n
\in \N)$ of elements of $\N^d$ such that $L^{(n)}_u=n $ for all
$n\geq 1$, and $L^{(n)}_v \rightarrow \infty$ as $n\rightarrow
\infty$, $\forall \,$ $v \in \{1, 2, \ldots, d\}\setminus\{u\}$, one has that
$\mu \left( L^{(n)} \right)\Rightarrow \nu.$ Then $\mu(L)\Rightarrow
\nu$ as $\min\{L_1, L_2, \ldots, L_d\}\rightarrow\infty$.
\end{proposition}
Let us specify that the proof of this proposition follows
exactly the proof given in \citet{Bradley3}, A2906 Proposition
(parts (I) and (III)) with just a small, insignificant
change.\\

For each $n\geq 1$ and each $\lambda \in [-\pi, \pi]$, define now the Fej{\'e}r kernel, $K_{n-1}(\lambda)$  by:
\begin{equation}\label{k11} K_{n-1}(\lambda):= \frac{1}{n} \left| \sum_{j=0}^{n-1} e^{ij\lambda}  \right|^2=\frac{\sin^2(n\lambda/2)}{n\sin^2(\lambda/2)}. \end{equation}
Elements of $[-\pi, \pi]^d$ will be denoted by $\vec{\lambda}:=(\lambda_1, \lambda_2, \ldots, \lambda_d)$.
For each $L \in \N^d$ define the ``multivariate Fej{\'e}r kernel" $G_{L}:[-\pi, \pi]^d\rightarrow [0,\infty)$ by:
\begin{equation}\label{k22} G_{L}(\vec{\lambda}):=\prod_{u=1}^{d} K_{L_u-1}(\lambda_u). \end{equation}

Also, on the ``cube" $[-\pi, \pi]^d$, let $m$ denote ``normalized Lebesque measure",
$m:=\text{Lebesque measure}/(2\pi)^d$.
\begin{lemma}\label{apt2}
Suppose $d$ is a positive integer. Suppose $f:[-\pi,\pi]^d\rightarrow \C$ is a continuous function. Then
\begin{equation*}\int_{\vec{\lambda} \in [-\pi,\pi]^d} G_{L}(\vec{\lambda}) \cdot f(\vec{\lambda}) dm(\vec{\lambda}) \rightarrow f(\vec{0}) \text{ as } \min\{L_1, L_2, \ldots, L_d\} \rightarrow\infty.\end{equation*}
\end{lemma}
Let us mention that Lemma \ref{apt2} is a special case of the multivariate Fej{\'e}r theorem, where the function $f$ is a periodic function with period $2\pi$ in every coordinate. For a proof of the one dimensional case, see \citet{Rud}, Theorem 8.15.\\

Further notations will be introduced and used throughout the entire paper.\\
If $a_n \in (0, \infty)$ and $b_n \in (0, \infty)$ for all $n\in \N$ sufficiently large, the notation $a_n\ll b_n$ means that $\limsup_{n\rightarrow\infty} a_n/b_n<\infty.$\\
If $a_n \in (0, \infty)$ and $b_n \in (0, \infty)$ for all $n\in \N$ sufficiently large, the notation $a_n\lesssim b_n$ means that $\limsup_{n\rightarrow\infty} a_n/b_n\leq 1.$\\
If $a_n \in (0, \infty)$ and $b_n \in (0, \infty)$ for all $n\in \N$ sufficiently large, the notation $a_n\sim b_n$ means that $\lim_{n\rightarrow\infty} a_n/b_n=1.$

\section{Central Limit Theorems } In this section we introduce two limit theorems that help us build up the main result, presented also in this section, as Theorem \ref{hilb}.

\begin{theorem}\label{r} Suppose $d$ is a positive integer. Suppose also that $X:=\left(X_k, k\in \Z^d \right)$
is a strictly stationary $\rho'$-mixing random field with the random
variables $X_k$ being real-valued such that $EX_0=0$ and
$EX^2_{0}<\infty$.

Then the following two statements hold:

(I)  The quantity
\begin{equation*} \sigma^2:=\lim_{\min\{L_1, L_2, \ldots, L_d\}\rightarrow\infty}\frac{ES^2(X, L)}{L_1\cdot L_2\cdot \ldots
\cdot L_d} \text{ exists in  }[0, \infty), \text{ and}
\end{equation*}
 (II) $\text{As } \min\{L_1, L_2, \ldots, L_d\}\rightarrow\infty,$
$(L_1\cdot L_2\cdot\ldots \cdot L_d)^{-1/2}S(X,L)\Rightarrow
N(0, \sigma^2).$ (Here and throughout the paper $\Rightarrow$ denotes convergence in distribution.)
\end{theorem}

\bpf The proof of the theorem has resemblance to arguments in earlier papers involving the $\rho^*$-mixing condition and similar properties as Theorem \ref{r} (see \citet{Brad92}
and \citet{Mil94}).
The proof will be written out for the case $d\geq 2$ since it is
essentially the same for the case $d=1$, but the notations for the
general case $d\geq 2$ are more complicated.

\textbf{Proof of (I).}  Our task is to show that there exists a number
$\sigma^2 \in [0, \infty)$ such that \beq\label{s2.2}
\lim_{\min\{L_1, L_2, \ldots, L_d\}\rightarrow\infty}\frac{ES^2
\left (X, L \right )}{L_1\cdot L_2\cdot \ldots \cdot L_d}=\sigma^2.
\eeq
For a given strictly stationary random field $X:=\left( X_k, k\in \Z^d \right)$ with mean zero and finite second moments, if $\rho'(n)\rightarrow 0$ as $n\rightarrow\infty$ then
$\zeta(n)\rightarrow 0$ as $n\rightarrow\infty$. Hence, by \cite{Bradley3} (Remark 29.4(V)(ii) and Remark 28.11(iii)(iv)), the random field $X$ has exactly one continuous spectral density function, $\sigma^2:=f(1, 1, \ldots, 1)$, where $f:[-\pi, \pi]^d\rightarrow [0,
\infty)$, and in addition, it is periodic with period $2\pi$ in every
coordinate.
In the following, by basic computations we compute the quantity given in \eqref{s2.2}. First we obtain that:
\begin{align}\label{s2.3}
\bsp & E  \left  |S  \left (X, L \right )   \right |^2=E \left
|\sum^{L_1}_{k_1=1}
\ldots\sum^{L_d}_{k_d=1}X_{(k_1,
\ldots, k_d)} \right |^2\\
& =\left (\sum^{L_1}_{k_1=1}
\ldots\sum^{L_d}_{k_d=1} \right)
\left(\sum^{L_1}_{l_1=1}
\ldots\sum^{L_d}_{l_d=1}\right )E
X_{(k_1, \ldots, k_d)}\overline{X_{(l_1, 
\ldots, l_d)}}.\\
\esp
\end{align}
We substitute the last term in the right-hand side of \eqref{s2.3} by the following expression (see \cite{Bradley3}, Section 0.19):
\begin{align}\label{cont2}
\bsp
& \frac{1}{(2\pi)^d} \left(\sum^{L_1}_{k_1=1}
\ldots\sum^{L_d}_{k_d=1} \right) \left
(\sum^{L_1}_{l_1=1}
\ldots\sum^{L_d}_{l_d=1}\right)
\int^\pi_{\lambda_1=-\pi}
\ldots
\int^\pi_{\lambda_d=-\pi}\\
&  e^{i  \left ((k_1-l_1)\lambda_1+
\ldots +
(k_d-l_d)\lambda_d  \right )}f(e^{i\lambda_1}, 
\ldots,
e^{i\lambda_d})d\lambda_d \ldots 
d\lambda_1  \\
& =\frac{1}{(2\pi)^d} \int^\pi_{\lambda_1=-\pi}
\ldots \int^\pi_{\lambda_d=-\pi}
f(e^{i\lambda_1},
\ldots, e^{i\lambda_d})\cdot\\
&\cdot \left
(\sum^{L_1}_{k_1=1}\sum^{L_1}_{l_1=1}e^{i(k_1-l_1)\lambda_1} 
\ldots
\sum^{L_d}_{k_d=1}\sum^{L_d}_{l_d=1}e^{i(k_d-l_d)\lambda_d} \right)
d\lambda_d \ldots  
d\lambda_1.
\esp
\end{align}
By \eqref{k11}, the right-hand side of \eqref{cont2} becomes:
\begin{align}\label{cont1}
\bsp
&\frac{1}{(2\pi)^d} \int^\pi_{\lambda_1=-\pi} 
\ldots
\int^\pi_{\lambda_d=-\pi} f(e^{i\lambda_1}, 
\ldots, e^{i\lambda_d})
\cdot\\
 & \cdot\frac{\sin^2 \left( L_1\lambda_1/2
 \right) }{\sin^2(\lambda_1/2)}
  \cdot \ldots \cdot
\frac{\sin^2  \left( L_d\lambda_d/2 \right) }{\sin^2(\lambda_d/2)}
d\lambda_d \ldots d\lambda_1 \\
&=\frac{1}{(2\pi)^d} \int^\pi_{\lambda_1=-\pi} 
\ldots
\int^\pi_{\lambda_d=-\pi} f(e^{i\lambda_1},
\ldots, e^{i\lambda_d})
\cdot\\
&\cdot(L_1 \cdot \ldots \cdot L_d) \cdot G_L(\lambda_1, \ldots, \lambda_d)d\lambda_d \ldots 
d\lambda_1,
\esp
\end{align}
therefore, by \eqref{s2.3}, \eqref{cont1} and the application of Lemma \ref{apt2}, we obtain that
\begin{align*}\label{s2.7} \bsp &\lim_{\min\{L_1, \ldots, L_d\}\rightarrow \infty}
\frac{ES^2 \left (X, \L \right )}{L_1 \cdot 
\ldots \cdot
L_d}=\lim_{\min\{L_1,\ldots, L_d\}\rightarrow \infty} \frac{1}{(2\pi)^d} \int^\pi_{\lambda_1=-\pi} \ldots
\int^\pi_{\lambda_d=-\pi}G_L(\lambda_1,\ldots, \lambda_d)\\
& \hspace{6cm} \cdot f(e^{i\lambda_1},
\ldots, e^{i\lambda_d}) d\lambda_d \ldots d\lambda_1\\
& \hspace{4.25cm}=f(1, \ldots, 1).
\esp \end{align*}
Hence, we can conclude that there exists a number $\sigma^2:=f(1, \ldots, 1)$
in $[0, \infty)$ satisfying \eqref{s2.2}. This completes the proof of part (I).

\textbf{Proof of (II).} Refer now to Proposition \ref{app1} from Section \ref{s:pre}.
Let $u \in \{1, 2, \ldots, d \}$ be arbitrary but fixed. Let $L^{(1)}, L^{(2)},
L^{(3)}, \ldots$ be an arbitrary fixed sequence of elements of
$\N^d$ such that for each $n\geq 1$, $L^{(n)}_u=n$ and $L^{(n)}_v
\rightarrow \infty$ as $n\rightarrow \infty$, $\forall \,$ $v \in \{1,
2, \ldots, d\}\setminus\{u\}$.
It suffices to show that \beq\label{g} \frac{S  \left(
X,L^{(n)} \right)}{\sqrt{L^{(n)}_1\cdot L^{(n)}_2\cdot \ldots \cdot
L^{(n)}_d}}  \Rightarrow N(0, \sigma^2) \text{ as }
n\rightarrow\infty. \eeq With no loss of generality, we can permute
the indices in the coordinate system of $\Z^d$, in order to have
$u=1$, and as a consequence, we have:
\begin{align}\label{defLn}
L^{(n)}_1=n \text{ for } n\geq1, \text{ and } L^{(n)}_v \rightarrow
\infty \text{ as } n\rightarrow \infty, \ \forall \, \ v \in \{2, \ldots, d\}.
\end{align}
Thus for each $n\geq1$, let us represent $L^{(n)}:=\left(
n, L^{(n)}_2, L^{(n)}_3, \ldots, L^{(n)}_d \right).$
We assume from now on, throughout the rest
of the proof that $\sigma^2>0$. The case $\sigma^2=0$ holds trivially by an application of Chebyshev Inequality.

\textbf{Step 1.} A common technique used in proving central limit theorems for random fields satisfying strong mixing conditions is the truncation argument whose effect makes the partial sum of the bounded random variables converge weakly to a normal distribution while the tails are negligible. To achieve this, for each integer $n\geq 1$, define the (finite) positive number
\beq\label{s2.15} c_n:= \left( L^{(n)}_2\cdot L^{(n)}_3 \cdot \ldots
\cdot
L^{(n)}_d \right)^{1/4}. \eeq
\eqref{defLn},
\beq\label{s2.16} c_n\rightarrow\infty \text{ as }
n\rightarrow\infty. \eeq
For each $n\geq 1$, we define the strictly stationary random field of bounded variables $X^{(n)}:= \left( X^{(n)}_k, k
\in \Z^d \right)$ as follows:
\beq\label{s2.17} \forall \, k \in \Z^d, \text{ } X^{(n)}_k
:=X_kI(|X_k|\leq c_n)-EX_0I(|X_0|\leq c_n). \eeq
Hence, by simple computations we obtain that $\forall \,
n\geq1$, \beq\label{s2.18} EX^{(n)}_0=0 \text{ and } Var X^{(n)}_0=E  \left (X^{(n)}_0  \right )^2
\leq EX^2_0 <\infty. \eeq
We easily also obtain that $\forall \, n\geq1$,
\beq\label{s2.21}
\left |X^{(n)}_0
  \right  |\leq 2c_n \text{ and } \left \|X^{(n)}_0 \right
\|_2\leq \|X_0  \|_2. \eeq
Next for $n\geq 1$, we define the strictly stationary random field of the tails of the $X_k$'s, $k \in \Z^d$, $\widetilde{X}^{(n)}:=
 \left (\widetilde{X}^{(n)}_k, k \in \Z^d  \right )$  
as follows (recall \eqref{s2.17} and the assumption $EX_0=0)$:
\beq\label{s2.44} \forall \, k \in \Z^d, \
\widetilde{X}^{(n)}_k :=X_k-X^{(n)}_k=X_kI(|X_k|> c_n)-EX_0I(|X_0|>
c_n). \eeq
As in \eqref{s2.44}, we similarly obtain by the dominated convergence theorem that \beq\label{s2.45} \forall \, n\geq1, \text{ }
E\widetilde{X}^{(n)}_0=0 \text{ and } E  \left (\widetilde{X}^{(n)}_0
\right )^2
\rightarrow 0 \text{ as } n \rightarrow\infty.\eeq
Note that $S\left(X, L^{(n)} \right):=\sum_{k}{X_k}=\sum_{k}{X^{(n)}_k}+\sum_{k}{\widetilde{X}^{(n)}_k}$, where all the sums are taken
over all $d$-tuples $k:=(k_1, k_2, \ldots, k_d) \in \N^d$ such that $1 \leq k_u \leq L_u \text{ for all } u\in \{1, 2, \ldots, d \}$.
Also, throughout the paper, unless specified, the notation $\sum_{k}$ will mean that the sum is taken over the same  set of indices as above.

\textbf{Step 2} (Parameters).
For each $ n \geq 1$, define the positive integer $q_n:=[n^{1/4}]$,
the greatest integer $\leq n^{1/4}$. Then it follows that \beq\label{s2.25} q_n\rightarrow\infty \text{ as }
n\rightarrow\infty. \eeq
Recall that $\rho'(X, n)\rightarrow 0$ as $n\rightarrow \infty$. As a consequence, we have
the following two properties: \beq\label{s2.22} \alpha(X,
n)\rightarrow 0 \text{ as } n\rightarrow \infty, \text{ and also}\eeq
\beq\label{s2.23} \text{there exists a positive integer j such that
} \rho'(X, j)< 1. \eeq
Let such a $j$ henceforth be fixed for the rest of the proof.
By \eqref{s2.22} and
\eqref{s2.25}, \beq\label{s2.27} \alpha(X, q_n)\rightarrow 0  \text{
as } n \rightarrow\infty. \eeq With $[x]$ denoting the greatest
integer $\leq x$, define the positive integers $m_n, \ n \geq
1$ as follows: \beq\label{s2.28} m_n:=\left [\min \left\{  q_n,
n^{1/10}, \alpha ^{-1/5}(X, q_n) \right\} \right ] .\eeq By the
equations \eqref{s2.28}, \eqref{s2.25}, and \eqref{s2.27}, we obtain
the following properties: \beq\label{s2.29} m_n  \rightarrow\infty
\text{ as } n \rightarrow\infty,  \eeq \beq\label{s2.30} m_n  \leq
q_n \text{ for all  }  n\geq 1,  \eeq
\beq\label{s2.32} \frac{m_n q_n}{ n}  \rightarrow 0 \text{ as } n
\rightarrow\infty, \text{ and}\eeq \beq\label{s2.33} m_n\alpha(X,
q_n)\rightarrow 0  \text{ as } n \rightarrow\infty. \eeq For each $
n \geq 1$, let $p_n$ be the integer such that \beq\label{s2.34}
m_n(p_n -1 +q_n) < n \leq m_n(p_n+q_n). \eeq
Hence we also have that \beq\label{s2-33} p_n\rightarrow\infty \text{ as
} n \rightarrow\infty \text{ and } m_np_n\sim n.\eeq
\textbf{Step 3} (The "Blocks"). In the following we decompose the partial sum of the bounded random variables $X^{(n)}_k, \ k \in \Z^d$ into ``big blocks" separated in between by ``small blocks". The ``lengths'' of both the big blocks and the small blocks, $p_n$ and $q_n$ respectively, have to
``blow up'' much faster than the (equal) numbers of big and
small blocks, $m_n$ (in addition to the fact that the ``lengths of the
``big blocks'' need to ``blow up'' much faster than the ``lengths''
of the ``small blocks''). This explains the way the positive integers $m_n$, $n\geq 1$ were defined in \eqref{s2.28}.
Referring to the definition of the random variables $X^{(n)}_k$ in
\eqref{s2.17}, for any $n \geq1$ and any two positive integers
$v\leq w $, define the random variable \beq\label{s2.36} Y^{(n)}(v,
w):= \sum _k  X^{(n)}_k ,\eeq \noindent where the sum is taken over
all $ k:=(k_1, k_2, \ldots, k_d) \in \N^d $  such that $v \leq k_1
\leq w$ and  $1 \leq k_u \leq L^{(n)}_u$ for all $ u\in \{2, \ldots, d \}$. Notice that for each $n\geq 1$, $S  \left (X^{(n)}, L^{(n)}
 \right )=Y^{(n)}(1, n).$
Referring to \eqref{s2.36}, for each $n\geq 1$, define the
random variables $U^{(n)}_k$ and $V^{(n)}_k$, as follows:
\begin{align}\label{s2.37}\bsp & \forall \, k \in \{ 1, 2, \ldots, m_n\}, \
U^{(n)}_k:= Y^{(n)} \left ((k-1)(p_n+q_n)+1, kp_n+(k-1)q_n
\right);\\
&  \text{(``big blocks")}\esp
\end{align}
    \begin{align}\label{s2.38}\bsp & \forall \, k \in \{ 1, 2, \ldots, m_n-1\}, \ V^{(n)}_k:= Y^{(n)}(kp_n+(k-1)q_n+1,
k(p_n+q_n));   \esp \end{align}\text{("small blocks")}, and \beq\label{s2.39}
V^{(n)}_{m_n}:= Y^{(n)}(m_np_n+(m_n-1)q_n+1, n). \eeq Note that by
\eqref{s2.30} and the first inequality in \eqref{s2.34}, for $ n\geq
1$,
\[ m_np_n+(m_n-1)q_n+1\leq
m_np_n+m_nq_n-m_n+1 \leq n. \] By \eqref{s2.36}, \eqref{s2.37},
\eqref{s2.38}, and \eqref{s2.39}, \beq\label{s2.40} \forall \, n\geq 1,
\text{ } S  \left (X^{(n)}, L^{(n)}  \right )=\sum^{m_n}_{k=1}
U^{(n)}_k + \sum^{m_n}_{k=1} V^{(n)}_k. \eeq
\textbf{Step 4} (Negligibility of the "small blocks").
Note that by \eqref{s2.38} and \eqref{s2.39},
$\sum^{m_n}_{k=1} V^{(n)}_k$ is the sum of  at most $m_n
\cdot q_n\cdot L^{(n)}_2 \cdot \ldots \cdot L^{(n)}_d$ of the random
variables $X^{(n)}_k$.
Therefore, by \eqref{s2.23} and  \cite{Bradley3},
Theorem 28.10(I), for any $n\geq 1$, the following holds:
\beq\label{s2.41} E \left |
\sum^{m_n}_{k=1} V^{(n)}_k  \right |^2 \leq C \left (m_n \cdot q_n
\cdot L^{(n)}_2 \cdot \ldots \cdot L^{(n)}_d  \right ) E  \left
(X^{(n)}_0  \right )^2,  \eeq where $C:=j^d \left (1+ \rho'(X, j) \right
)^d/\left (1- \rho'(X, j) \right )^d, $
and as a consequence,  by \eqref{s2.32} and \eqref{s2.18},
we obtain that
\begin{align}\label{s2.42}
\bsp & E \left | \frac{  \sum^{m_n}_{k=1} V^{(n)}_k
}{\sigma   \sqrt{n \cdot L^{(n)}_2 \cdot \ldots \cdot L^{(n)}_d}
 }\right |^2
 \leq \frac{C (m_n  q_n)E  \left (X^{(n)}_0  \right
)^2}{n\cdot\sigma^2} \rightarrow 0 \text{ as } n\rightarrow \infty.\esp
\end{align}
Hence, the ``small blocks" are negligible:  \begin{equation}\label{neg}\frac{\sum^{m_n}_{k=1} V^{(n)}_k}{\sigma \sqrt{n \cdot
L^{(n)}_2 \cdot \ldots \cdot L^{(n)}_d}}\rightarrow 0 \text{ in
probability as } n \rightarrow\infty.\end{equation}
By an obvious analog of \eqref{s2.42}, followed by \eqref{s2.45},
for each $n\geq 1$, we obtain that \begin{align}\label{cca}\frac{\sum_{k} \widetilde{X}^{(n)}_k}{ \sigma \sqrt{n \cdot
L^{(n)}_2 \cdot \ldots \cdot L^{(n)}_d}}\rightarrow 0 \text{ in
probability as } n \rightarrow\infty.\end{align}
\textbf{Step 5} (Application of the Lyapounov CLT).
For a given $n\geq 1$, by the definition of $U^{(n)}_k$ in \eqref{s2.37} and the
strict stationarity of the random field $X^{(n)}$, the random
variables $U^{(n)}_1, U^{(n)}_2, \ldots, U^{(n)}_{m_n} $ are
identically distributed.
For each $n\geq 1$, let $\widetilde{U}^{(n)}_1,
\widetilde{U}^{(n)}_2, \ldots, \widetilde{U}^{(n)}_{m_n} $ be
independent, identically distributed random variables whose common
distribution is the same as that of $U^{(n)}_1$. Hence,
since  $\forall \, n\geq1,
 EX^{(n)}_0=0$, we have the following:
\[E \widetilde{U}^{(n)}_1=E U^{(n)}_1=0 \text{ and }  Var  \left (\sum^{m_n}_{k=1} \widetilde{U}^{(n)}_k
\right ) =m_n E  \left ( \widetilde{U}^{(n)}_1  \right )^2 = m_n E \left ( U^{(n)}_1  \right )^2.\]
By \eqref{s2.23}, we can refer to \cite{Bradley3}, Theorem 29.30, a result which gives us a Rosenthal inequality for $\rho'$-mixing random fields. Also, using the fact that $EU^2_1 \sim \sigma^2\left(
p_n\cdot L^{(n)}_2 \cdot \ldots \cdot L^{(n)}_d  \right)$ (see \eqref{s2.2}), together with the equations \eqref{s2.21},  \eqref{s2.18},
and assuming without loss of
generality that $EX^2_0\leq 1$, the following holds:
\begin{align}\label{s2.58}
\bsp & \frac{E \left (U^{(n)}_1 \right )^4}{m_n \left (EU^2_1 \right
)^2} \lesssim
\frac{ C_R \left (p_n \cdot L^{(n)}_2 \cdot \ldots
\cdot L^{(n)}_d \cdot E   \left |X^{(n)}_0   \right |^4
+ \left (p_n \cdot L^{(n)}_2 \cdot
\ldots \cdot L^{(n)}_d \cdot EX^2_0  \right)^2  \right ) }{m_n p^2_n \sigma^4\left ( L^{(n)}_2 \cdot \ldots \cdot
L^{(n)}_d \right)^2 } \\
&\hspace{2cm} \leq \frac{16C_R  p_n c^4_n\left ( L^{(n)}_2 \cdot \ldots \cdot
L^{(n)}_d \right )} {m_n p^2_n \left ( L^{(n)}_2 \cdot \ldots \cdot
L^{(n)}_d \right)^2 \sigma^4} + \frac{C_R p^2_n \left ( L^{(n)}_2
\cdot \ldots \cdot L^{(n)}_d \right )^2  }{m_n p^2_n \left (
L^{(n)}_2 \cdot \ldots \cdot L^{(n)}_d \right )^2
\sigma^4} \\
& \hspace{2cm} \leq \frac{16C_R } {m_n p_n  \sigma^4} + \frac{C_R
}{m_n
\sigma^4}\rightarrow 0
\text{ as }  n \rightarrow\infty \text{ by \eqref{s2-33} and \eqref{s2.29}}.
\esp
\end{align}
Since $U_1-U^{(n)}_1$ is the sum of
$p_n\cdot L^{(n)}_2 \cdot \ldots \cdot L^{(n)}_d$ random variables
$\widetilde{X}^{(n)}_k$, applying an obvious analog of \eqref{s2.41}, followed by \eqref{s2.2} and \eqref{s2.45}, we have that as $n\rightarrow\infty$,
\begin{equation*}
\label{s2.60}
\frac{E \left (U_1-U^{(n)}_1 \right )^2}{EU^2_1} \lesssim
\frac{Cp_n \left (L^{(n)}_2 \cdot \ldots \cdot L^{(n)}_d \right )E
\left (\widetilde{X}^{(n)}_0 \right )^2}{p_n \left ( L^{(n)}_2 \cdot
\ldots \cdot L^{(n)}_d \right )\sigma^2}=\frac{C  E \left (\widetilde{X}^{(n)}_0 \right
)^2}{\sigma^2} \rightarrow 0.\end{equation*}
As a consequence, after an application of
Minkowski Inequality to the quantity $ \left |\|U_1\|_2- \left \| U^{(n)}_1  \right \|_2  \right  |/ \|U_1\|_2$, we have that
\beq\label{s2.61}   E \left (U^{(n)}_1 \right)^2  \sim
EU^2_1.
\eeq
Hence, by \eqref{s2.58} and \eqref{s2.61}, the following holds:
\begin{align*}\label{s2.62}
\bsp & \frac{E \left (U^{(n)}_1 \right )^4}{m_n \left (E \left(
U^{(n)}_1  \right)^2 \right )^2} \sim \frac{E\left (U^{(n)}_1 \right
)^4}{m_n(EU^2_1)^2}
\rightarrow 0 \text{ as } n\rightarrow\infty.
 \esp
\end{align*}
Therefore, due to Lyapounov CLT  (see
\citet{Bill2}, Theorem 27.3), it follows that  \begin{equation}\label{nou}
\left (\sqrt{m_n} \left \|U^{(n)}_1  \right \|_2 \right )^{-1}
\sum^{m_n}_{k=1} \widetilde{U}^{(n)}_k \Rightarrow N(0, 1) \text{ as } n \rightarrow \infty.
 \end{equation}
\textbf{Step 6.} As in \citet{Bradley3}, Theorem 29.32,
we similarly obtain by \eqref{s2.36}, \eqref{s2.37} and \eqref{s2.33} that as $n \rightarrow\infty$,
\[ \sum ^{m_n-1}_{k=1} \alpha \left (\sigma   \left (U^{(n)}_j, 1\leq j \leq k  \right ),
\sigma \left (U^{(n)}_{k+1} \right ) \right ) \leq \sum
^{m_n-1}_{k=1} \alpha \left (X^{(n)},q_n \right ) \leq m_n \alpha(X,
q_n)\rightarrow 0.
\]
Hence, by \eqref{nou} and by \cite{Bradley3}, Theorem 25.56, the following
holds: \beq\label{s2.64} \left (\sum ^{m_n}_{k=1} U^{(n)}_k \right )
/ \left (\sqrt{m_n}  \left \|U^{(n)}_1   \right \|_2 \right )
\Rightarrow N(0, 1)  \text{ as } n \rightarrow\infty . \eeq
Refer to the first sentence of Step 5. For each $n\geq 1$,
 \beq\label{s2.68}
E \left (\sum ^{m_n}_{k=1} U^{(n)}_k  \right )^2= m_n E \left
(U^{(n)}_1 \right )^2 + 2\sum ^{m_n-1}_{k=1}\sum ^{m_n}_{j=k+1}
EU^{(n)}_k U^{(n)}_j.
 \eeq
Using similar arguments as in \citet{Bradley3}, Theorem 29.31 (Step 9), followed by \eqref{s2.58} and \eqref{s2.61}, and
\eqref{s2-33}, $E \left (U^{(n)}_1 \right )^4/\left (E \left(
U^{(n)}_1  \right)^2 \right )^2  \rightarrow C_R/\sigma^4$ as $n\rightarrow\infty$. Hence we obtain that
$\left  \|U^{(n)}_1   \right \|^2_4 \ll E
 \left (U_1^{(n)}  \right )^2.$
 As a consequence,
by \eqref{s2.68},
\begin{equation}\label{nnou}    \left \|\sum ^{m_n}_{k=1} U^{(n)}_k  \right \|^2  \sim \left( m_n E  \left (U^{(n)}_1  \right )^2\right)^{1/2}.\end{equation}
Applying an obvious analog of \eqref{s2.41} for
$S \left (\widetilde{X}^{(n)},
L^{(n)} \right ):= S \left (X, L^{(n)} \right )- S \left (X^{(n)},
L^{(n)} \right )$, followed by \eqref{s2.2} and \eqref{s2.45}, the following holds:
\begin{align}\label{nouu}
\bsp & E \left (S \left (\widetilde{X}^{(n)}, L^{(n)} \right )
\right )^2 / E \left (S \left (X, L^{(n)} \right ) \right )^2
 \lesssim
C  E \left (\widetilde{X}^{(n)}_0 \right )^2 / \sigma^2 \rightarrow 0 \text{ as } n \rightarrow\infty. \esp
\end{align}
Using Minkowski Inequality for  $\left |  \left \|S \left (X, L^{(n)} \right ) \right \|_2 -
\left \| S \left (X^{(n)}, L^{(n)} \right ) \right \|_2
\right | / \left \|S \left (X, L^{(n)} \right ) \right \|_2$, by \eqref{nouu} it follows that
\begin{equation}\label{nnou1} \left \| S \left (X^{(n)}, L^{(n)} \right ) \right \|_2
 \sim\left
 \|S \left (X, L^{(n)} \right ) \right \|_2.\end{equation}
Now apply again Minkowski Inequality for \[\left | \left \|  \sum^{m_n}_{k=1} U^{(n)}_k   \right \|_2 - \left
\| S \left (X^{(n)}, L^{(n)} \right ) \right \|_2 \right | / \left
\| S \left (X^{(n)}, L^{(n)} \right ) \right \|_2,  \] and by the formulation of $S \left (X^{(n)}, L^{(n)} \right )$ given in \eqref{s2.40}, followed by \eqref{s2.41}, \eqref{nnou}, \eqref{s2.2} and by \eqref{s2.32}, we obtain that
\begin{equation}\label{s2.67}
\left \| S \left (X^{(n)}, L^{(n)} \right ) \right \|_2 \sim  \left
\| \sum^{m_n}_{k=1} U^{(n)}_k \right \|_2. \end{equation}
Hence, by \eqref{nnou} and \eqref{nnou1},
\begin{equation*}\label{sum1}\left\| S
\left( X,L^{(n)} \right) \right\|_2 \sim \left( m_n E  \left (U^{(n)}_1  \right )^2\right)^{1/2}. \end{equation*}
As a consequence, by \eqref{s2.64} and the fact that $\left\| S
\left( X,L^{(n)} \right) \right\|_2 \sim \sigma \sqrt{ n \cdot L^{(n)}_2
\cdot \ldots \cdot L^{(n)}_d }$ (see \eqref{s2.2}), it follows the following:
\begin{equation}\label{nor}\frac{ \sum_{k=1}^{m_n} U_k^{(n)}}{\sigma \sqrt{ n \cdot L^{(n)}_2
\cdot \ldots \cdot L^{(n)}_d }} \Rightarrow N(0, 1) \text{ as } n
\rightarrow\infty.\end{equation}
{\bf Step 7.} Refer to the definition of $S \left (X^{(n)}, L^{(n)} \right )$ given in \eqref{s2.40}. By \eqref{neg} and \eqref{nor}, followed by \cite{Bradley3}, Theorem 0.6, we obtain the following weak convergence:
\begin{equation}\label{nor2}\frac{ S \left (X^{(n)}, L^{(n)} \right ) }{\sigma \sqrt{n \cdot
L^{(n)}_2 \cdot \ldots \cdot L^{(n)}_d }} \Rightarrow N(0, 1) \text{
as } n \rightarrow\infty.\end{equation}
Refer now to the definition of $S \left (X, L^{(n)} \right )$ given just after \eqref{s2.45}. By another application of Theorem 0.6 from \cite{Bradley3} for \eqref{cca} and \eqref{nor2}, we obtain that \eqref{g} holds, and hence, the proof of (II) is complete. Moreover, the proof of the
theorem is complete.
\epf
\begin{theorem}\label{rm} Suppose $d$ and $m$ are each a positive integer. Suppose $X:=(X_k, k\in
\Z^d)$ is a strictly stationary $\rho'$-mixing random field  with
$X_k:=(X_{k1}, X_{k2}, \ldots, X_{km})$ being (for each $k$) an $m$-dimensional random vector
such that $\forall \, i \in \{ 1, 2 ,
 \cdots, m\}$, $X_{ki}$ is a real-valued random variable
with $EX_{ki}=0$ and $EX^2_{ki}<\infty$.

Then the following statements hold:\\
(I) For any $i \in \{ 1, 2 , \ldots, m\}$, the quantity
\begin{equation*} \sigma_{ii}= \lim_{min\{L_1, L_2, \ldots, L_d\}\rightarrow\infty}
\frac{ES^2_{L,i}}{L_1\cdot L_2 \cdot\ldots\cdot L_d}  \text{ exists
in  }[0, \infty), \text{}
\end{equation*}
where for each $L \in \N^d$ and each $i \in \{1, 2, \ldots, m\}$,
\begin{align}\label{s3.5}
\bsp &S_{L,i}:=\sum_{k}{X_{ki}},  \esp
\end{align}
with the sum being
taken over all $k:=(k_1, k_2, \ldots, k_d) \in \N^d$ such that
$1 \leq k_u \leq L_u $ for all  $u\in \{1, 2, \ldots, d \}$.\\
(II) Also, for any two distinct elements $i$, $j$ $\in \{ 1, 2 , \ldots, m\}$,
\begin{equation*} \gamma(i, j)= \lim_{min\{L_1, L_2, \ldots, L_d\}\rightarrow\infty} \frac{E(S_{L,i}- S_{L,j})^2}{L_1\cdot L_2
\cdot\ldots\cdot L_d}\text{ exists in  }[0, \infty).
\end{equation*}
(III) $\text{ Furthermore, as } \min\{L_1, L_2, \ldots, L_d\}\rightarrow\infty,$
\begin{equation*} \frac{S(X, L)}{\sqrt{L_1\cdot
L_2\cdot \ldots \cdot L_d}}\Rightarrow N(0_m, \Sigma),\text{ where}
\end{equation*}\noindent
\beq\label{s3.2}
 \Sigma:=(\sigma_{ij}, 1 \leq i \leq j \leq m) \text{ is the } m\times m
\text{ covariance matrix defined by } \eeq
\beq\label{s3.3} \text{for } i\neq j, \text{ } \sigma_{ij}=
\frac{1}{2}(\sigma_{ii}+ \sigma_{jj}- \gamma(i, j)), \eeq with
$\sigma_{ii}$ and $\gamma(i, j)$ defined in part (I), respectively
in part (II).\\
 (The fact that the matrix $\Sigma$ in (III) is symmetric and nonnegative definite (and can therefore be a covariance matrix), is part of the conclusion of (III).)
\end{theorem}
\bpf A distant resemblance to this theorem is a bivariate central limit theorem of \citet{Mil95}. The proof of Theorem \ref{rm} will be divided in the following
parts:\\
\textbf{Proof of (I) and (II). } Since $\sigma_{ii}$, respectively
$\gamma(i, j)$ exist by Theorem \ref{r}(I), parts (I) and (II) hold.\\
\textbf{Proof of (III). } For the clarity of the proof, the strategy used to prove this part is the following:\\
(i) It will be shown that the matrix $\Sigma$ defined in part (III) is symmetric and nonnegative definite.\\
(ii) One will then let $Y:=(Y_1, Y_2, \ldots, Y_m)$ be
a centered normal random vector with covariance matrix $\Sigma$, and the task will be to show that \beq\label{s3.6}
\frac{S(X, L)}{\sqrt{L_1\cdot L_2 \cdot\ldots\cdot
L_d}}\Rightarrow Y \text{ as } \min\{L_1, L_2, \ldots, L_d\}\rightarrow\infty. \eeq
(iii) To accomplish that, by the Cramer-Wold Device Theorem (see
\citet{Bill2}, Theorem 29.4) it suffices to show that for an arbitrary $t \in \R^m$,
 \beq\label{s3.7} t \cdot
\frac{S_{L}}{\sqrt{L_1\cdot L_2 \cdot\ldots\cdot
L_d}}\Rightarrow t \cdot Y \text{ as } \min\{L_1, L_2, \ldots, L_d\}\rightarrow\infty, \eeq where ``$\cdot$" denotes the scalar product.

Let us first show (i). In order to achieve this task, let us introduce $\Sigma^{(L)}:=\left( \sigma^{(L)}_{ij}, 1 \leq i
\leq j \leq m \right)$ to be the $m\times m$ covariance matrix
defined by
\beq\label{s3.9} \sigma_{ij}^{(L)}= ES_{L,i}S_{L,j}=\frac{1}{2}
\left( ES^2_{L,i}+ES^2_{L,j}-E(S_{L,i}-S_{L,j} \right)^2). \eeq
Note that $\sigma_{ii}^{(L)}= ES_{L,i}^2$ for $i \in \{1, 2, \ldots, m\}$.
Our main goal is to prove that
\beq\label{s3.15} \lim_{\min\{L_1, L_2, \ldots, L_d\}\rightarrow\infty} \frac{\Sigma^{(L)}}{L_1\cdot L_2
\cdot\ldots\cdot L_d} = \Sigma \text{ (defined in } \eqref{s3.2}
).\eeq It actually suffices to show that
\beq\label{s3.16} \lim_{\min\{L_1, L_2, \ldots, L_d\}\rightarrow\infty} \frac{\sigma_{ij}^{(L)}}{L_1\cdot L_2
\cdot\ldots\cdot L_d}
 = \sigma_{ij} \text{, } \forall \, \text{ } 1 \leq i \leq j
\leq m.
 \eeq
By the definition of $\sigma_{ij}^{(L)}$ given in \eqref{s3.9}, followed by the distribution of the limit (each of the limits exist by Theorem \ref{rm}, parts (I) and (II)), the
left-hand side of \eqref{s3.16} becomes:
\begin{align*}\label{s3.17}
\bsp & \frac{1}{2} \lim_{\min\{L_1, L_2, \ldots, L_d\}\rightarrow\infty} \frac{1}{L_1\cdot L_2 \cdot\ldots\cdot L_d}
\left( ES^2_{L,i}+ES^2_{L,j}-E \left( S_{L,i}-S_{L,j} \right)^2 \right)  \\
&=\frac{1}{2} \left( \sigma_{ii}+ \sigma_{jj}- \gamma(i, j) \right) =
\sigma_{ij}.
 \esp
\end{align*}
Let us recall that each of these limits exist by Theorem \ref{rm}, parts (I) and (II).
Hence, \eqref{s3.16} holds. As a consequence, \eqref{s3.15} also
holds.

In the following, one should mention
that since $\Sigma^{(L)}$ is the $m \times m$ covariance matrix of
$S_{L,i}$, one has that $\Sigma^{(L)}$ is symmetric and nonnegative
definite. That is, $\forall \, r:=(r_1, r_2, \ldots, r_m) \in \R^m,
\ r \Sigma^{(L)} r^{'} \geq 0 $.
Therefore, $\forall \, r \in \R^m,$  $r ( L_1\cdot L_2
\cdot\ldots\cdot L_d)^{-1} \Sigma^{(L)} r^{'} \geq 0, $ and
moreover,
\[ \forall \, r \in \R^m, \text{ } r \left( \lim_{\min\{L_1, L_2, \ldots, L_d\}\rightarrow\infty} (L_1\cdot L_2 \cdot\ldots\cdot L_d)^{-1}
\Sigma^{(L)} \right) r^{'} \geq 0.\] By \eqref{s3.15}, we get that
$ \forall \, r \in \R^m, \text{ } r \Sigma r^{'} \geq 0$, and hence,
$\Sigma$ is also symmetric (trivially by \eqref{s3.15}) and nonnegative definite. Hence, there exists a centered normal random vector $Y:=(Y_1,
Y_2, \ldots, Y_m)$ whose covariance
matrix is $\Sigma$, and therefore, the proof of (i) is complete.

(ii) Let us now take $Y:=(Y_1, Y_2, \ldots, Y_m)$ be
a centered normal random vector with covariance matrix $\Sigma$, defined in \eqref{s3.2}. As we mentioned above, the task now is to show that \eqref{s3.6}
holds. In order to accomplish this task, by part (iii), one would need to show \eqref{s3.7}.

(iii) So, let $t:=(t_1, t_2, \ldots, t_m)$ be an arbitrary fixed element of $\R^m$.
 We can notice now that
\begin{align}\label{s3.8}
\bsp & t \cdot S_{L}=\sum_{i=1}^m t_i S_{L,i}, \text{ where }
S_{L,i} \text{ is defined in } \eqref{s3.5}.\esp
\end{align}
  We can also notice that
$t\cdot X_1, t\cdot X_2, \ldots $ is a strictly
stationary $\rho'$-mixing random sequence with real-valued random
variables that satisfy $E \left(
t\cdot X_1 \right)=t \cdot E X_1=t \cdot 0_m=0$,
and $ E \left( t\cdot X_1 \right)^2 <
\infty$.
For these random variables we can apply Theorem \ref{r}.
Therefore, we
obtain that as $ \min\{L_1, L_2, \ldots, L_d\}\rightarrow\infty$,   \beq\label{s3.12}  t \cdot
\frac{S_{L}}{\sqrt{L_1\cdot L_2 \cdot\ldots\cdot L_d}}
\Rightarrow N(0, \sigma^2), \eeq where \beq\label{s3.13}
\sigma^2:=\lim_{\min\{L_1, L_2, \ldots, L_d\}\rightarrow\infty}
\frac{ E \left( t \cdot S_{L}
 \right)^2}{L_1\cdot L_2 \cdot\ldots\cdot L_d}. \eeq
Moreover, by \eqref{s3.8}, \eqref{s3.9}, and \eqref{s3.15}, \eqref{s3.13} becomes:
\begin{align}\label{s3.14}
\bsp &\sigma^2=\lim_{\min\{L_1, L_2, \ldots, L_d\}\rightarrow\infty}
\frac{ E \left( \sum_{i=1}^m t_i S_{L,i} \right)^2}{L_1\cdot L_2
\cdot\ldots\cdot L_d} \\
& \, \ = \lim_{\min\{L_1, L_2, \ldots, L_d\}\rightarrow\infty}
\frac{1}{L_1\cdot L_2 \cdot\ldots\cdot L_d} \left( \sum_{i=1}^m
t_i^2
ES_{L,i}^2 +  \right.\\
&  \left. +\sum_{1 \leq i < j \leq m} t_i t_j
   \left( ES^2_{L,i}+ES^2_{L,j}-E \left( S_{L,i}-S_{L,j} \right)^2  \right)  \right)  \\
&=  t  \left( \lim_{\min\{L_1, L_2, \ldots, L_d\}\rightarrow\infty}
\frac{\Sigma^{(L)}}{L_1\cdot L_2 \cdot\ldots\cdot L_d}  \right)
t^{'}=t\Sigma t^{'}.
 \esp
\end{align}
By \eqref{s3.12} and \eqref{s3.14}, one can
conclude that \beq\label{s3.20} t \cdot
\frac{S_{L}}{\sqrt{L_1\cdot L_2 \cdot\ldots\cdot L_d}}
\Rightarrow N \left( 0, t \Sigma t^{'} \right) \text{ as }
\min\{L_1, L_2, \ldots, L_d\}\rightarrow\infty. \eeq
Also, since the random vector $Y$ is centered normal with
covariance matrix $\Sigma$, one has that $t \cdot Y$ is a
normal random variable with mean 0 and variance ($1\times 1$
covariance matrix) $t \Sigma t^{'}$.
Hence, by \eqref{s3.20}, \eqref{s3.7} holds, therefore \eqref{s3.6}
holds.  This completes the proof of Theorem \ref{rm}.
\epf
\begin{theorem}\label{hilb}
Suppose $H$ is a separable real Hilbert space, with inner product
$\langle \cdot    ,\cdot \rangle$ and norm $\|\cdot \|_H$. Suppose
$X:=(X_k, k\in \Z^d)$ is a strictly stationary $\rho'$-mixing random
field with the random variables $X_k$ being H-valued, such that
\beq\label{s4.6} EX_{{0}}=0_H \text{ and} \eeq
\beq\label{s4.7} E \left\| X_{{0}}  \right\|_H^2
<\infty. \eeq
Suppose $\{ e_i \}_{i\geq 1}$ is an orthonormal basis of $H$ and that $X_{ki}:= \langle X_k ,e_i
\rangle$ for each pair $(k, i).$

Then the following statements hold:\\
(I) For each $i \in \N$, the quantity
\begin{equation*} \sigma_{ii}= \lim_{min\{L_1, L_2, \ldots, L_d\}\rightarrow\infty}
\frac{ES^2_{L,i}}{L_1\cdot L_2 \cdot\ldots\cdot L_d}  \text{ exists
in  }[0, \infty), \text{ where}
\end{equation*}
\begin{equation}\label{s4.11}
S_{L,i}:=\sum_{k}{X_{ki}}, \text{ the sum being
taken over all } k:=(k_1, k_2, \ldots, k_d) \in \N^d \end{equation}
such that
 $1 \leq k_u \leq L_u$  for all $u\in \{1, 2, \ldots, d \}$.\\
(II) Also, for any two distinct elements, $i$, $j$ $\in \N$,
\begin{equation*} \gamma(i, j)= \lim_{\min\{L_1, L_2, \ldots, L_d\}\rightarrow\infty} \frac{E(S_{L,i}- S_{L,j})^2}{L_1\cdot L_2
\cdot\ldots\cdot L_d}\text{ exists in  }[0, \infty).
\end{equation*}
(III) $\text{ Furthermore, as } \min\{L_1, L_2, \ldots, L_d\}\rightarrow\infty,$
\begin{equation*} \frac{S(X, L)}{\sqrt{L_1\cdot
L_2\cdot \ldots \cdot L_d}}\Rightarrow N \left( 0_H,
\Sigma^{(\infty)} \right),\end{equation*}
where the ``covariance operator" $\Sigma^{(\infty)}:=(\sigma_{ij}, i\geq 1,  j\geq 1)$ is symmetric, nonnegative definite, has finite trace and it is defined by
\beq\label{s4.9} \text{for } i\neq j, \ \sigma_{ij}=
\frac{1}{2}(\sigma_{ii}+ \sigma_{jj}- \gamma(i, j)), \eeq with
$\sigma_{ii}$ and $\gamma(i, j)$ defined in part (I), respectively
in part (II). (Recall that $\Rightarrow$ denotes convergence in distribution and also the statement before Lemma \ref{hl1}.)
\end{theorem}
\bpf
The proof of the theorem will be divided in the following parts:\\
\textbf{Proof of (I) and (II). } Since $\sigma_{ii}$, respectively
$\gamma(i, j)$ exist by Theorem \ref{r}(I), parts (I) and (II) hold.\\
\textbf{Proof of (III). } The rest of the proof will be divided into
five short steps, as follows:\\
\textbf{Step 1.} Since the Hilbert space $H$ is separable, one can
consider working with  the  separable Hilbert space $l_2$.
 Let us recall that $\forall \, k\in \Z^d$, $X_k=(X_{k1}, X_{k2}, X_{k3}, \cdots)$ is an
$l_2$-valued random variable with real-valued components such that
\beq\label{s4.12} EX_{ki}=0, \ \forall \, i\geq 1 \text{ and}
\eeq
\beq\label{s4.13} E\|X_k \|_H^2 <\infty. \eeq
For any given $m \in \N$, if one considers the first $m$ coordinates of the $l_2$-valued
random variable $X_k$, $X^{(m)}_k:=(X_{k1}, X_{k2}, \ldots, X_{km})$, by Theorem \ref{rm} we obtain:
\beq\label{s4.14} \frac{S_{L}^{(m)}}{\sqrt{L_1\cdot L_2
\cdot\ldots\cdot L_d}}\Rightarrow N \left( 0_m, \Sigma^{(m)}
\right) \text{ as } \min\{L_1, L_2, \ldots, L_d\}\rightarrow\infty,
\eeq \noindent where $\Sigma^{(m)}:=(\sigma_{ij}, 1 \leq i \leq j
\leq m)$ is the
$m\times m$ covariance matrix defined as in \eqref{s3.2}. Let us specify that here and below, for any given $L \in \N^d$ and $m\in \N$, the random variable $S^{(m)}_L$ is defined by:
\begin{equation*}S_L^{(m)}:=\sum_{k}{X_{k}^{(m)}}, \text{ the sum being
taken over all } k:=(k_1, k_2, \ldots, k_d) \in \N^d\end{equation*}
such that
 $1 \leq k_u \leq L_u$  for all $u\in \{1, 2, \ldots, d \}$.

\textbf{Step 2.} Suppose $m \in \N$. Let $\widetilde{Y}^{(m)}:= \left( Y_1^{(m)},
Y_2^{(m)}, \ldots , Y_m^{(m)}  \right)$ be an $\R^{m}$-valued random
vector whose distribution on $(\R^m, \mathcal{R}^m) $ is $N \left(
0_m, \Sigma^{(m)} \right), \Sigma^{(m)} $ being the same covariance
matrix defined in \eqref{s3.2}. By Step 1, we have that \beq\label{s4.15} \frac{S_{L}^{(m)}}{\sqrt{L_1\cdot L_2
\cdot\ldots\cdot L_d}}\Rightarrow  \widetilde{Y}^{(m)} \text{ as }
\min\{L_1, L_2, \ldots, L_d\}\rightarrow\infty. \eeq
Let $\mu _m$ be the probability measure on $(\R^m, \mathcal{R}^m) $
of the random vector $\widetilde{Y}^{(m)}$ and let $\mu _{m+1}$ be
the probability measure on $(\R^{m+1}, \mathcal{R}^{m+1}) $ of the
random vector $\widetilde{Y}^{(m+1)}:= \left( Y_1^{(m+1)},
Y_2^{(m+1)}, \ldots , Y_m^{(m+1)}, Y_{m+1}^{(m+1)}  \right)$, whose
distribution is $N \left( 0_{m+1}, \Sigma^{(m+1)}  \right)$.
One should specify that $ \Sigma^{(m+1)}:=(\sigma_{ij}, 1 \leq i
\leq j \leq m+1) $ is the $(m+1)\times (m+1)$  covariance matrix
defined in \eqref{s3.2}, where the integer $m$ in \eqref{s3.2} corresponds to $m+1$ here.
\begin{Cl}\label{hcl1} For each $m \in \N$, $ \left( Y_1^{(m+1)}, Y_2^{(m+1)}, \ldots ,
Y_m^{(m+1)} \right)$ (that is, the first $m$ coordinates of the random vector $\widetilde{Y}^{(m+1)}$) has the same distribution as
$\widetilde{Y}^{(m)}:= \left( Y_1^{(m)}, Y_2^{(m)}, \ldots ,
Y_m^{(m)}  \right). $
\end{Cl}
\bpf Since the random vector $\tilde{Y}^{(m+1)}$ is (multivariate) centered normal, it follows automatically that $ \left( Y_1^{(m+1)}, Y_2^{(m+1)}, \ldots ,
Y_m^{(m+1)} \right)$ (the first $m$ coordinates) is centered normal.
For the two centered normal random vectors $\widetilde{Y}^{(m)}$ and see above $ \left( Y_1^{(m+1)}, Y_2^{(m+1)}, \ldots ,
Y_m^{(m+1)} \right)$, the $m\times m$ covariance matrices are the same (with the common entries being the elements $\sigma_{ii}$ and $\sigma_{ij}$ defined in
Theorem \ref{rm}). From this observation, as well as the fact that a (multivariate) centered normal distribution is uniquely determined by its covariance matrix, Claim \ref{hcl1} follows.
\epf
Now, by Kolmogorov's Existence Theorem (see \citet{Bill2}, Theorem
36.2), there exists on some probability space $(\Omega,
\mathcal{F}, P)$ a sequence of random variables $Y:=(Y_1, Y_2, Y_3,
\ldots )$ such that for each $m\geq 1$, the $m$-dimensional random vector $(Y_1,
Y_2, \ldots, Y_m)$ has distribution $\mu _m$  on $(\R^m,
\mathcal{R}^m) $.
\begin{Cl}\label{hcl2} $Y$ is a centered normal $l_2$-valued random variable.\end{Cl}
\bpf First of all, one should prove that $Y$ is an $l_2$-valued
random variable, whose (random) norm has a finite second moment; that is, \beq\label{s4.16} E\|Y \|^2_{l_2} <\infty.
\eeq More precisely, one should check that \beq\label{s4.17}
\sum_{i=1}^{\infty}EY_i^2= \sum_{i=1}^{\infty} \sigma_{ii} <\infty ,
\text{ where } \sigma_{ii} = Cov(Y_i, Y_i)= EY_i^2.   \eeq
Since for every $i\geq 1$, $S_{L,i}$ is the sum of $L_1\cdot L_2
\cdot\ldots\cdot L_d$ real-valued random variables $X_{ki}$, by
an obvious analog of \eqref{s2.41}, followed by the definition of $\sigma_{ii}$, given in
part (I) of the theorem, we obtain the following inequality:
\beq\label{s4.19} \sigma_{ii}\leq C\cdot E|X_{0i}|^2, \text{ where }
C \text{ is the constant defined just after }\eqref{s2.41}   \eeq (with $j\geq 1$ fixed such that $\rho'(X, j)<1).$
Therefore, by \eqref{s4.19} and \eqref{s4.13},
 \[ \sum_{i=1}^{\infty} \sigma_{ii} \leq C\sum_{i=1}^{\infty}
E|X_{0i}|^2< \infty .\] Hence, \eqref{s4.17} holds, that is $Y$ is
an $l_2$-valued random variable, whose (random) norm has a finite second moment.
In order to prove that $Y$ is a {\em normal} $l_2$-valued random variable,
it now suffices to show the following:
\begin{align}\label{s4.20} \bsp &\forall \, m\geq 1  \text{ and } \forall \, (r_1,
r_2, \ldots, r_m) \in \R^m, \text{ the real-valued random variable}\\
&\sum_{i=1}^{m} r_iY_i \text{ is normal (possibly degenerate)}. \esp
\end{align}
In order to show \eqref{s4.20}, let $m\geq 1$ and $(r_1, r_2, \ldots, r_m) \in \R^m$. As we mentioned earlier, for each $m\geq 1$,
the random vector $(Y_1, Y_2, \ldots, Y_m)$ is centered normal with
covariance matrix $\Sigma^{(m)}$, defined in \eqref{s3.2}.
Therefore, $\sum_{i=1}^{m} r_iY_i$ is a centered normal real random
variable.
Hence, $Y$ is a centered normal $l_2$-valued random variable (possibly
degenerate) whose ``covariance operator" is defined in 
\eqref{s4.9}, and therefore, the proof of Claim \ref{hcl2} is complete.
 \epf
\textbf{Step 3.} Refer now to Proposition \ref{app1} from Section \ref{s:pre}.
Let $u \in \{1, 2, \ldots, d \}$ be arbitrary but fixed. Let $L^{(1)}, L^{(2)},
L^{(3)}, \ldots$ be an arbitrary fixed sequence of elements of
$\N^d$ such that for each $n\geq 1$, $L^{(n)}_u=n$ and $L^{(n)}_v
\rightarrow \infty$ as $n\rightarrow \infty$, $\forall \,$ $v \in \{1,
2, \ldots, d\}\setminus\{u\}$.

Suppose $m\geq1$.
Consider the following sequence: \begin{equation*}\frac{S^{(m)}
\left( X, L^{(1)} \right) }{\sqrt{L^{(1)}_{1}\cdot L^{(1)}_{2}
\cdot\ldots\cdot L^{(1)}_{d}}}, \frac{S^{(m)} \left( X,
L^{(2)}\right) }{\sqrt{L^{(2)}_{1}\cdot L^{(2)}_{2} \cdot\ldots\cdot
L^{(2)}_{d}}}, \ldots,  \frac{S^{(m)} \left( X,L^{(n)}\right)
}{\sqrt{L^{(n)}_{1}\cdot L^{(n)}_{2} \cdot\ldots\cdot L^{(n)}_{d}}},
\ldots.\end{equation*} By Step 1, one has the following:
\beq\label{s4.21} \frac{S^{(m)} \left( X,L^{(n)}\right)
}{\sqrt{L^{(n)}_{1}\cdot L^{(n)}_{2} \cdot\ldots\cdot
L^{(n)}_{d}}}\Rightarrow N \left( 0_m, \Sigma^{(m)}  \right) \text{
as } n\rightarrow \infty,\eeq \noindent where $\Sigma^{(m)}$ is the
$m \times m$ covariance matrix defined in \eqref{s3.2}.

\textbf{Step 4.} Let $\mathcal{P}$ denote the family of distributions of the $l_2$-valued random variables $S_{L}/\sqrt{L_1\cdot L_2 \cdot\ldots\cdot L_d}$, $L \in \N^d$.
By Lemma \ref{hl1}, in order to show that $\mathcal{P}$ is tight, one
should show that
\begin{align}\label{s4.22} \bsp &\lim_{N\rightarrow
\infty} \sup_{L \in \N^d} E  \left( \sum_{i=N}^{\infty} \left \langle
\frac{S_{L}}{\sqrt{L_1\cdot L_2 \cdot\ldots\cdot L_d}} ,e_i \right
\rangle ^2 \right)=0,\esp
\end{align}
as well as the fact that for $N=1$ the supremum in \eqref{s4.22} is finite.

Let $N\geq1$ and $L \in \N^d$. Then using \eqref{s4.11}, followed by an obvious analog of \eqref{s2.41}, we obtain the following:
\begin{align*}\label{s4.23} \bsp & E  \left( \sum_{i=N}^{\infty} \left
\langle \frac{S_{L}}{\sqrt{L_1\cdot L_2 \cdot\ldots\cdot L_d}} ,e_i
\right \rangle ^2 \right)
= \frac{1}{L_1\cdot L_2 \cdot\ldots\cdot L_d}
\sum_{i=N}^{\infty} ES_{L,i} ^2
 \leq C \sum_{i=N}^{\infty}
E|X_{0i}|^2.
\esp
\end{align*}
Since $E \|X_{0}\|_{H}^2<\infty$,
one has that \begin{equation}\label{lim0}
\lim_{N \rightarrow\infty }\sum_{i=N}^{\infty}E|X_{0i}|^2=0.\end{equation} Also by \eqref{s4.7}, for $N=1$ the sum in \eqref{lim0} is finite.
Hence \eqref{s4.22} holds, and as a consequence,
$\mathcal{P}$ is tight.
Moreover, $\mathcal{P}$ is tight along the sequence $L^{(1)},
L^{(2)}, L^{(3)}, \cdots$, hence the family of distributions $
 \left\{ S \left( X,L^{(n)}  \right) /\sqrt{L^{(n)}_1\cdot L^{(n)}_2 \cdot\ldots\cdot
L^{(n)}_d}  \right\}$ is tight. As a consequence, the sequence $S \left(
X,L^{(n)}  \right)/\sqrt{L^{(n)}_1\cdot L^{(n)}_2 \cdot\ldots\cdot
L^{(n)}_d}$ contains a weakly convergent subsequence.

\textbf{Step 5. } Let $Q$ be an infinite set in $ \N$. Assume that as $n\rightarrow\infty$, $n\in Q$,
the sequence
$S \left( X,L^{(n)}
\right)/\sqrt{L^{(n)}_1\cdot L^{(n)}_2 \cdot\ldots\cdot L^{(n)}_d}
\Rightarrow W:=(W_1, W_2, W_3, \ldots)$.

By Step 3,
$(W_1, W_2, \ldots, W_m)$ is $N \left( 0_m, \Sigma^{(m)}  \right)$,
where $\Sigma^{(m)}:=(\sigma_{ij}, 1 \leq i \leq j \leq m)$ is the
$m\times m$ covariance matrix defined in \eqref{s3.2}.
Hence, the distribution of the random vector $(W_1, W_2, \ldots, W_m)$ is the same as the distribution of $Y^{(m)}, \forall \,
m$. Thus the distributions of $W$ and $Y$ are identical.
Therefore,
\begin{equation}\label{great} \frac{S \left( X,L^{(n)}
\right)}{\sqrt{L^{(n)}_1\cdot L^{(n)}_2 \cdot\ldots\cdot L^{(n)}_d}}
\Rightarrow Y \text{ as } n\rightarrow\infty, \text{ } n\in Q.
\end{equation}
Hence, we obtain that the convergence in \eqref{great} holds along the entire sequence of positive integers, and as a consequence,
\begin{equation*}\frac{S(X, L)}{\sqrt{L_1\cdot L_2 \cdot\ldots\cdot L_d}} \Rightarrow Y
\text{ as } \min\{L_1, L_2, \ldots, L_d\}.\end{equation*} Therefore, part (III)
holds, and hence, the proof of the theorem is complete. \epf

{\bf{Acknowledgment}}

The result is part of the author's Ph.D. thesis at Indiana University (Bloomington). The author thanks her advisor, Professor Richard Bradley, to whom she is greatly indebted for his advice and support not only in this work but also during the graduate years at Indiana University.

\bibliographystyle{alea2}

\begin{thebibliography}{3}
\bibitem[Billingsley(1995)]{Bill2}P. Billingsley. \emph{ Probability and Measure}, 3rd ed., Wiley, New-York (1995).
\bibitem[Billingsley(1999)]{Bill1}P. Billingsley. \emph{Convergence of Probability Measures}, 2nd ed., Wiley, New-York (1999).
\bibitem[Bradley(2007)]{Bradley3}R. C. Bradley. \emph{ Introduction to Strong Mixing Conditions}, volumes 1, 2, 3, Kendrick Press, Heber City, Utah (2007).
\bibitem[Bradley(1992)]{Brad92}R. C. Bradley. On the spectral density and asymptotic normality of weakly dependent random fields, \emph{J. Theor. Probab.} 5, 355-373 (1992).
\bibitem[Dedecker and Merlev\`{e}de(2002)]{DM} J. Dedecker and F. Merlev\`{e}de. Necessary and suffcient conditions for the conditional
central limit theorem, \emph{Ann. Probab.} 30 1044-1081(2002).
\bibitem[Laha and Rohatgi(1979)]{Laha}R. G. Laha and V. K. Rohatgi. \emph{ Probability Theory}, Wiley (1979).
\bibitem[Mal'tsev and  Ostrovskii(1982)]{Mal}V. V. Mal'tsev and E. I. Ostrovskii. Central limit theorem for strictly stationary processes in Hilbert space, \emph{Theor. Probab. Appl.} 27, 357-359 (1982).
\bibitem[Merlev\`{e}de(2003)]{FM}F. Merlev\`{e}de. On the central limit theorem and its weak invariance principle for strongly mixing sequences with values in a Hilbert space via martingale approximation, \emph{J. Theor. Probab.} 16, 625-653 (2003).
\bibitem[Merlev\`{e}de, Peligrad, and Utev(1997)]{MPU}F. Merlev\`{e}de, M. Peligrad, and  S. Utev. Sharp conditions for the central limit theorem of linear processes in a Hilbert space, \emph{J. Theor. Probab.} 10, 681-693 (1997).
\bibitem[Miller(1994)]{Mil94}C. Miller. Three theorems on $\rho^*$-mixing random fields, \emph{J. Theor. Probab.} 7, 867-882 (1994).
\bibitem[Miller(1995)]{Mil95}C. Miller. A central limit theorem for the periodograms of a $\rho^*$-mixing random field, \emph{Stochastic Process. Appl.} 60, 313-330 (1995).
\bibitem[Rudin(1974)]{Rud}W. Rudin. \emph{ Principles of Mathematical Analysis}, 2nd ed., McGraw-Hill, New York (1974).
\end{thebibliography}

\end{document}